\documentclass[a4paper,11pt]{amsart}

\usepackage{enumerate}
\usepackage{parskip}
\usepackage{amsmath,amssymb,colordvi}
\usepackage{color}
\usepackage{graphics}
\usepackage{graphicx}
\newtheorem{theorem}{Theorem}[section]

\newtheorem{lemma}[theorem]{Lemma}

\newcommand{\di}{\displaystyle}
\newcommand{\R}{{\Bbb R}}

\newcommand{\N}{{\Bbb N}}
\newcommand{\C}{{\Bbb C}}
\definecolor{B}{cmyk}{1,1,0,0}
\definecolor{R}{rgb}{0.8,0.0,0.0}

 \numberwithin{equation}{section}

\begin{document}

\title[]{Exponential dichotomy for a class of asymptotically autonomous delayed differential equations}

\author{Hel\'i Elorreaga}
\address{University of B\'io-B\'io, Department of Mathematics}
\email{helorreaga@ubiobio.cl}
\thanks{}

\author{Adri\'an G\'omez}
\address{University of B\'io-B\'io, Department of Mathematics}
\curraddr{}
\email{agomez@ubiobio.cl}
\thanks{}



\maketitle

\begin{abstract}
In this work we give a criterion to have an exponential dichotomy over all $\mathbb{R}$ for delayed systems $x'(t)=L(t)x_t$, where $L_{\pm}=\lim_{t\to\pm\infty}L(t)$,  and the systems $x'(t)=L_{\pm}x_t$ are autonomous and hyperbolic. The proof is based in a suitable choice of weighted Sobolev spaces $L^p(\R,\mu)$ and  $W^{1,p}(\R,\mu)$, where the involves differential operators become Fredholm's operators and this property are strongly exploited.

\end{abstract}

\section{Introduction and main result}
Our main goal in this paper is to prove the existence of an exponential dichotomy over all $\mathbb{R}$,  for a class of asymptotically autonomous delayed differential equation, with $N$ discrete delays. This type of dichotomies play a role in the study of stability of linear non-autonomous systems and also in some nonlinear systems which are linearized around a solution. The stability  is  an important topic in the investigation of systems of differential equa\-tions, as it is well known, in the case  of a hyperbolic autonomous  system, it  is determined by
the analysis of the roots of its associated characteristic equation. Nevertheless, there are some non-autonomous systems where their sta\-bility  cannot be determined by looking at the roots of their associated characteristic equation. Indeed, L. Markus and H.
Yamabe  in \cite{Markus-Yamabe} constructed a counterintuitive example. Therefore another approach is needed for the investigation of the stability of non-autonomous linear system. Coppel in \cite{Coppel}, using the concept of ordinary and exponential dichotomy, which generalizes the hyperbolicity in the non-autonomous case,   establishes new results in stability theory, becoming a fruitful tool among those alternative approaches. In this way, is important to mention the works  Palmer \cite{PalmerHG, Palmer, Palmer2} for ordinary differential equation and Lin \cite{LIN} for functional differential equations. We also refer to the work of Sacker and Sell \cite{Sacker} in which was studied properties related with exponential dichotomies of a general class of infinite-dimensional systems.\\
 As is explain in \cite{Barreira-valls-hyperbolicity}, there is an extensive literature on the relation between the concept of \emph{admissibility} and stability.    Recently, in \cite{Barreira-valls-hyperbolicity}  Luis Barreira and Claudia Valls   stablished a caracterization of exponential dichotomies via the concept of admisibility of suitable space. After in \cite{Barreira-valls-perturbations}, using as hypothesis the existence  of an exponential dichotomy on $\R$ for a nonautonomous linear homogeneous delay equation, they showed a functional version of the famous Hartman-Grobman Theorem, in the sense of Palmer \cite{PalmerHG}.\\
More specifically, let $C$ the Banach space of all continuous functions over $[-r,0]$, provided with the supremum norm: $C=C([-r,0],\R^n)$, with the delay $r>0$. Let the homogeneous delayed equation

\begin{equation}\label{eq0}
x'(t)=L(t)x_t,\ t\in\R,
\end{equation}
with $x_t(\theta):=x(t+\theta)$, for $\theta\in[-r,0]$ and $x_t\in C$ and  $L(t):C\to \R^n$ are con\-tinuous linear operators for all $t\in\R$. \\ 
In \cite{Barreira-valls-hyperbolicity} is showed that the existence of a exponential dichotomy over all $\R$ can be caracterized by the admisibility of the spaces $(C_b,M)$, where 
\begin{eqnarray*}
C_b&:=& \{f:\mathbb{R}\to\mathbb{R}^n / f \mbox{ is bounded and continuous}\}, \mbox{and}  \\
  M&:=& \{f:\mathbb{R}\to\mathbb{R}^n/ f \mbox{ is measurable and }\int_t^{t+1}|f(s)|ds<\infty, \quad \mbox{ for all } t\in\mathbb{R}\}.
\end{eqnarray*}  
(See Theorem 2, \cite{Barreira-valls-hyperbolicity}.\\ However, in the practice is not easy to check the admissibility of $(C_b,M)$.
 \\
 Following this point of view, in this work we prove that   the admissibility of our spaces $(W^{1,p}(\mathbb{R}, d\mu),L^p(\mathbb{R}, d\mu))$, with a suitable  finite measure $\mu$ given by
 \begin{equation}\label{medida}
    \frac{d\mu}{dt} = \frac{1}{1+t^2},
\end{equation}
  is equivalent to  the existence of an exponential dichotomy for (\ref{eq0}) on $\R$, for a family of $L(t)$ asymptotically hyperbolic functionals.
Moreover,  defining the differential operator $\Lambda_L:W^{1,p}(\mathbb{R},d\mu)\to L^p(\mathbb{R},d\mu)$ by
\begin{equation}
(\Lambda_Lx)(t):=x'(t)-L(t)x_t=x'(t)-\sum\limits_{j=0}^{N}A_j(t)x(t-r_j),
\end{equation} we show that the admissibility of our pair of weighted spaces depends of the Fredholm's properties of $\Lambda_L$, which in their turn are contained in the following non-trivial adaptation of Theorem A in \cite{Mallet-Paret}:

\begin{theorem}\label{T-Fredholm}
Assuming that the homogeneous equation (\ref{eq0}) is asymptotically hyperbolic and taking the finite measure $\mu$. Then

\begin{enumerate}

\item The operator $\Lambda_L: W^{1,p}_{\mu}\to L_{\mu}^{p}$ is a Fredholm operator,  for any $1\leq p\leq \infty$,
\item  The kernels $\mathcal{K}_L \subseteq W^{1,p}_{\mu}$ and $\mathcal{K}_{L^*} \subseteq W^{1,q}_{\mu}$ of $\Lambda_L$ and $\Lambda_{L^*}$ respectively are independent of $p, q$, with $p^{-1}+q^{-1}=1$
\item The range $\mathcal{R}_L^p\subseteq L^p_{\mu}$ of $\Lambda_L$ in $L_{\mu}^p$ is given by  
\begin{equation}\label{range}
    \mathcal{R}_L^p = \Big\{h\in L_{\mu}^p : \int_{\mathbb{R}}\overline{y(t)} h(t)d\mu = 0, \mbox{ for all } y\in\mathcal{K}_{L^{*}}\Big\}.
\end{equation}
\end{enumerate}
In particular,
\begin{equation}\label{index}
    \dim \mathcal{K}_{L^{*}} = \textit{codim }  \mathcal{R}_L^p, \quad \dim \mathcal{K}_L = \textit{codim } \mathcal{R}_{L^{*}}^p, \quad \mbox{Ind}(\Lambda_L)=-\mbox{Ind}(\Lambda_{L^{*}}). 
\end{equation}
\end{theorem}

Supported in the last theorem, we show  the following characterization of the exponential dichotomy: 
\newpage 
\begin{theorem}\label{T-E-D}
Assuming that the homogeneous equation (\ref{eq0})  is asymptotically hyperbolic, and 
\begin{enumerate}
\item $\mbox{Ind}(\Lambda_L)=0$, 
\item $\Lambda_L: W^{1,p}_{\mu}\to L_{\mu}^{p}$ has a trivial kernel $ \mathcal{K}_L$  
\end{enumerate} 
Then, the  equation (\ref{eq0}) has an exponential dichotomy on $\mathbb{R}$.
\end{theorem}

This article is organized as follow: In the next section we summary some basic fact and the notation used throughout the work. In sections 3 we show Theorem (\ref{T-Fredholm}). For this  aim, we  obtain a Green's function  associated with the operator $\Lambda_L$ and estimation of it in Lemma (\ref{Prop-G(t,z)}), which in turns is  based on the autonomous case in Lemma  (\ref{T-Green-L_0}). Also  an important bounds over the solutions of (\ref{eq1}) and (\ref{eq2}) appear in Lemma (\ref{P-cotas-x}). Finally, in section 4, we show the Theorem (\ref{T-E-D}).


\section{Preliminaries and notation}

We start defining our weighted spaces 
\begin{equation}
    L^p(\mathbb{R},d\mu) = \{f:\mathbb{R}\to \R^n : f \mbox{ is measurable and } \int_{\mathbb{R}}|f(t)|^pd\mu<\infty\},\nonumber
\end{equation}
for  $1\leq p<\infty$, 
\begin{eqnarray}
    L^{\infty}(\mathbb{R},d\mu) = \{f:\mathbb{R}\to \R^n :& &  f \mbox{ is measurable and there is a constant $C$ such that }\nonumber\\
      & & |f(t)|\frac{d\mu}{dt}\leq C \mbox{ a.e. on } \mathbb{R}\},
\end{eqnarray}
and
\begin{equation}
    W^{1,p}(\mathbb{R},d\mu) = \{f\in L^p(\mathbb{R},d\mu) : f \mbox{ is absolutely continuous and } f'\in L^p(\mathbb{R},d\mu)\}.\nonumber
\end{equation}
For simplicity, we call $L^p_\mu$, $W^{1,p}_\mu$ to our weighted spaces and $L^p$, $W^{1,p}$ to the classic spaces without weight, over all $\R$.
\\

Writing the equation (\ref{eq0}) in the more familiar way
\begin{equation}\label{eq1}
    x'(t)=\sum\limits_{j=0}^{N}A_j(t)x(t-r_j),
\end{equation}
where, as we saw  
\begin{equation}
 r=r_N>r_{N-1}>\cdots >r_1>r_0=0    
\end{equation}
are the delays and the matrices $A_j(t)$ are measurable,  uniformly bounded on $\R$ with $\lim\limits_{t\to\pm\infty}=A_{\pm,j}\in M_n(\R)$. Associated with (\ref{eq1}) we consider the non-homogeneous equation

\begin{equation}\label{eq2}
    x'(t) = \sum\limits_{j=0}^{N}A_j(t)x(t-r_j)+h(t),
\end{equation}
where  $h:\mathbb{R}\to \mathbb{R}^n$ is a function in  $ L^p(\mathbb{R},d\mu)$.\\
A continuous function $x:J'\to \mathbb{R}^n$ defined on the larger interval
\begin{equation}
    J'=\{t+\theta : t\in J \mbox{ and } \theta\in[-r,0]\}\nonumber
\end{equation}
is a solution of (\ref{eq1}) on an interval $J$ if $x$ is absolutely continuous on $J$, and satisfies (\ref{eq1}) for almost every $t\in J$. Thus, for such $x$, we have that $x_t\in C$ for every $t\in J$. 


The linear operator $\Lambda_L: W^{1,p}_\mu\to L^p_\mu$  associated with the equation (\ref{eq1}) is defined by 
\begin{equation}\label{Lambda_L}
    (\Lambda_{L}x)(t)=x'(t) - \sum\limits_{j=0}^{N}A_j(t)x(t-r_j).
\end{equation}

Is easy to see that, since the coefficients $A_j$  have an uniform bound, the operator  $\Lambda_L$ is bounded for $1\leq p\leq\infty$. 

We recall that a bounded operator $A:X\to Y$, where $X, Y$ are Banach spaces is called a  \emph{Fredholm Operator} if $\dim(\ker(A))$ and $\dim(\ker(A'))$ are finite and  the range $\mathcal{R}(A)$ is a closed subespace of $Y$. If $A$ is  Fredholm operator, we define its index as\\ 
$$\mbox{Ind}(A)=\dim(\ker(A))-\dim(\ker(A'))$$ (see, for example \cite{ms}, pp 101).\\

\subsection*{Exponential Dichotomy}
Let $x$ be the unique solution of (\ref{eq1}) with initial data  $x_s=\phi$, for  $\phi\in C([-r, 0], \mathbb{R}^n)$. The \emph{solution operator} $T(t,s):C([-r, 0], \mathbb{R}^n)\to C([-r,0], \mathbb{R}^n)$ is given by
\begin{equation}\label{Solution-op}
    T(t,s)\phi=x_t, \ \ \mbox{ for } t\geq s \mbox{ and } \phi\in C.
\end{equation}
We will say that the equation (\ref{eq1})  has an \emph{Exponential Dichotomy} on an interval $J$, with positive constants $D$ and $\lambda$ if there are strongly continuous projections $P(s)$, $Q(s)=I-P(s)$, $s\in J$, satisfying
\begin{enumerate}
    \item[(i)] $T(t,s)P(s)=P(t)T(t,s)$, for $t\geq s$ in $J$;
    \item[(ii)] $\overline{T}(t,s)=T(t,s)|_{\mathcal{R}Q(s)}$ is an isomorphism of $\mathcal{R}Q(s)$ onto $\mathcal{R}Q(t)$ for $t\geq s$, and $\overline{T}(s,t):\mathcal{R}Q(t)\to\mathcal{R}Q(s)$ is defined as the inverse of $\overline{T}(t,s)$;
    \item[(iii)] $\|T(t,s)P(s)\|\leq De^{-\lambda(t-s)}$, for $t\geq s$ in $J$;
    \item[(iv)] $\|\overline{T}(s,t)Q(t)\|\leq De^{-\lambda(t-s)}$, for $t\geq s$ in $J$.
\end{enumerate}

\newpage
\subsection*{Adjoint Equation for (\ref{eq1})}

Calling   \begin{eqnarray*} \omega(t)&:=&(1+t^2)^{-1},\\
                             k(t)&:=&-\frac{2t}{1+t^2},\\ 
                             M^{\pm}_j(t)&:=&\frac{1+t^2}{1+(t\pm r_j)^2}
\end{eqnarray*}
 is easy to see that 
\begin{enumerate}
    \item[($P_1$)] $\omega'(t)=k(t)\omega(t)$.
    \item[($P_2$)] $w(t\pm r_j)= M^{\pm}_j(t) \omega(t)$. 
\end{enumerate}
With that in mind, we can calculate the adjoint of  $\Lambda_L$  in the following

\begin{lemma}
Let  $\Lambda_L$ the linear operator defined in (\ref{Lambda_L}). Its  adjoint ope\-ra\-tor $\Lambda^{*}_L$ is given by 
\begin{equation}\label{adjunta-Lambda_L}
    (\Lambda_L^{*} y)(t) = -y'(t) - \sum\limits_{j=0}^{N}B^{*}_j(t) y(t+r_j)
\end{equation}
for $y\in W^{1,q}_{\mu}$, where the matrices $B^{*}_j$ are the transposed conjugated of  
\begin{equation}\label{B_j}
  \begin{array}{rcl}
      B_0(t) &=& k(t)I + A_0(t)  \\
      B_j(t) &=& M^{+}_j(t) A_j(t+r_j), \qquad \mbox{ for } j=1,2,...,N 
  \end{array}  
\end{equation}
\end{lemma}
\begin{proof}
Let $x\in W^{1,p}_{\mu}$ and $y\in W^{1,q}_{\mu}$, where $1\leq p,q\leq \infty$ with $\displaystyle\frac{1}{p}+\frac{1}{q}=1$. We recall that the duality between these weighted spaces is given by 
$$\langle x,y \rangle_\mu:= \int_{\mathbb{R}}\overline{y(t)}x(t)d\mu=\int_{\mathbb{R}}\overline{y(t)}x(t)\omega(t)dt.$$

Using that,   integrating by parts and applying the  $(P_1)$ property,
\begin{eqnarray}
    \int_{\mathbb{R}} \overline{y(t)}x'(t)\omega(t)dt &=& -\int_{\mathbb{R}} \overline{y'(t)}\omega(t)x(t)dt - \int_{\mathbb{R}}\overline{y(t)}\omega'(t)x(t)dt\nonumber\\
    &=& -\int_{\mathbb{R}}\overline{y'(t)}x(t)\omega(t)dt - \int_{\mathbb{R}} \overline{y(t)}x(t)k(t)\omega(t)dt\nonumber\\    
&=& -\int_{\mathbb{R}}\overline{y'(t)}x(t)d\mu - \int_{\mathbb{R}}\overline{y(t)}x(t)k(t)d\mu .\label{IPP}
\end{eqnarray}
Now, we shall prove that
\begin{equation}\label{igualdad-adj}
    \int_{\mathbb{R}}\overline{y(t)}(\Lambda_L x)(t)d\mu = \int_{\mathbb{R}}\overline{(\Lambda_L^{*}y)(t)}x(t)d\mu
\end{equation}
From (\ref{IPP}) and $(P_2)$ :
\begin{eqnarray}
   \int_{\mathbb{R}}\overline{y(t)}(\Lambda_L x)(t)d\mu &=& \int_{\mathbb{R}}\overline{y(t)}x'(t)d\mu - \int_{\mathbb{R}}\overline{y(t)}\sum\limits_{j=0}^{N}A_j(t)x(t-r_j)d\mu \nonumber\\
   &=& -\int_{\mathbb{R}}\overline{y'(t)}x(t)d\mu - \int_{\mathbb{R}}\overline{y(t)}x(t)k(t)d\mu\nonumber\\
   && -\sum\limits_{j=0}^{N}\int_{\mathbb{R}}\overline{y(t)}A_j(t)x(t-r_j)w(t)dt\nonumber\\
   &=& -\int_{\mathbb{R}}\overline{y'(t)}x(t)d\mu - \int_{\mathbb{R}}\overline{y(t)}x(t)k(t)w(t)dt\nonumber \\
   && -\sum\limits_{j=0}^{N}\int_{\mathbb{R}}\overline{y(t+r_j)}A_j(t+r_j)x(t)w(t+r_j)dt\nonumber \\
   &=& -\int_{\mathbb{R}}\overline{y'(t)}x(t)d\mu - \int_{\mathbb{R}}\overline{y(t)}[k(t)I+A_0(t)]x(t)w(t)dt\nonumber\\
   && - \sum\limits_{j=1}^{N}\int_{\mathbb{R}}\overline{y(t+r_j)}A_j(t+r_j)x(t)M^{+}_{j}(t)w(t)dt
\end{eqnarray}
Defining the matrices $B_j$ as in (\ref{B_j}) we get
\begin{equation}
    \int_{\mathbb{R}}\overline{y(t)}(\Lambda_L x)(t)d\mu = \int_{\mathbb{R}}\Big[\overline{-y'(t)-\sum\limits_{j=0}^{N}B^{*}_j(t)y(t+r_j)}\Big]x(t)d\mu\nonumber
\end{equation}
as we desired.
\end{proof}

In view of the last Lemma, we define the adjoint equation of (\ref{eq1}) by
\begin{equation}\label{adjunta_of_L}
    y'(t)=L^{*}(t)y_t=-\sum\limits_{j=0}^{N}B^{*}_j(t)y(t+r_j)
\end{equation}
where
\begin{equation}
    L^{*}(t)\varphi=-\sum\limits_{j=0}^{N}B^{*}_j(t)\varphi(r_j), \qquad \varphi\in C([0,r],\mathbb{C}^d)\nonumber
\end{equation}
and the matrices $B_j$ defined as in (\ref{B_j}).
 Note that \begin{equation}
    (\Lambda_L^{*}y)(t)=-y'(t) + L^{*}(t)y_t\nonumber
\end{equation}
which implies that $$\Lambda_L^{*} = -\Lambda_{L^{*}}.$$

\subsection{Asymptotically Hyperbolic Systems}
The characteristic equation associa\-ted with the autonomous system

\begin{equation}\label{equautohomo}
x'(t)=L_0x_t=\sum\limits_{j=0}^{N}A_{j,0}x(t-r_j),
\end{equation}
 is given by $\det \Delta_{L_0}(s)=0$, where 
\begin{equation}\label{polcarac}
    \Delta_{L_0}(s)=sI-\sum_{j=0}^{N}A_{j,0}e^{-sr_j}.
\end{equation}
The system (\ref{equautohomo}) is \emph{hyperbolic} if
\begin{equation}\label{hip-cond}
    \det\Delta_{L_0}(iz) \neq 0, \ \ z\in\mathbb{R}.
\end{equation}

When we can write the non-autonomous operator $L(t)$  as a sum
\begin{equation}\label{L0+M(t)}
    L(t)= L_0 + M(t)
\end{equation}
of a constant coefficient operator defined in (\ref{equautohomo}) and a perturbation term $M(t)$ such that 

\begin{equation}
    \lim\limits_{t\to \pm\infty}\|M(t)\|=0,
\end{equation}
we say that (\ref{eq0}) is \emph{asymptotically autonomous} at $\pm\infty$ respectively. In terms of operators, we say that $\Lambda_L$ is asymptotically autonomous at $\pm\infty$ respectively if 
\begin{equation}
\Lambda_L=\Lambda_{L_0}+M(t),
\end{equation}
where, abusing the notation, is interpreted as operator $M(t):W^p_\mu\to L^p_\mu$  defined as $M(t)(x):=M(t)x_t$, with $M(t)$ the functional perturbation in (\ref{L0+M(t)}) and 

    $$\lim\limits_{t\to \pm\infty}\|M(t)\|=0,$$
in the suitable operator norm.  We can write this functional in the same way that $L(t)$ as
$$M(t)\phi:=\sum\limits_{j=0}^NC_j(t)\phi(-r_j),$$
for some measurable, uniformly bounded matrix functions $C_j(t)$ satisfying $$\lim_{t\to\pm \infty} C_j(t)= 0.$$
Of course we can mix these concepts and talk about a system \emph{asymptotically } \\ \emph{hyperbolic} in $\pm\infty$ respectively when it has both properties described before in the corresponding limit.\\
In the case of a system $x'(t)=L(t)x_t$ asymptotically hyperbolic in both limits $\pm\infty$ we say  that it is \emph{asymptotically hyperbolic}. Note that the limit systems  at $\pm\infty$ need not be the same.

\section{Proof of Theorem \ref{T-Fredholm}}

Supposing that the system (\ref{equautohomo}) is hyperbolic, in \cite{Mallet-Paret}, section 4  is used the inverse Fourier transformation
\begin{equation}\label{G0}
    G_0(t)=\frac{1}{2\pi}\int_{\R} e^{itz}\Delta^{-1}_{L_0}(iz)dz,
\end{equation}
in order to produce a Green's function for this system, we also need this tool, and in the following theorem, we show that it remains true in the case of our  weighted spaces. 
\newpage
\begin{lemma}\label{T-Green-L_0}
Assuming that the operator $L_0$ given by (\ref{equautohomo}) is  hyperbolic, we have 

\begin{enumerate}
\item The function $G_0(t)$ satisfies the following exponential estimation:
\begin{equation}\label{G0-estimate}
    |G_0(t)|\leq K_0e^{-a_0|t|}, 
\end{equation}
with $t\in\mathbb{R}$, and $K_0>0$ and $a_0>0$.
\item The associated operator  $\Lambda_{L_0}$ is an isomorphism from $W_{\mu}^{1,p}$ onto $L_{\mu}^p$,  for $1\leq p\leq\infty$.
\item The  inverse operator $\Lambda^{-1}_{L_0}$ is  given by the convolution
\begin{equation}
    (\Lambda^{-1}_{L_0}h)(t) = (G_0*h)(t)=\int_{\R} G_0(t-z)h(z)d\mu(z)\nonumber
\end{equation}

\end{enumerate} 
\end{lemma}

\begin{proof} For convenience of the reader, we divide the proof in several steps and  repeat and adapt some argument in  the proof of Theorem 4.1 \cite{Mallet-Paret} which still work in our case. 
\begin{enumerate}
\item[\textbf{Step 1}]\emph{$G_0(t)$ is an absolutely continuous function on $\R-\{0\}$, with a jump discontinuity satisfying $$G_0(0+)-G_0(0-)=I.$$}

In fact, considering $G_0(t)$ as a tempered distribution, and applying Fourier transform on  

\begin{equation}
G_0'(t)-\sum_{j=0}^NA_{j,0}G_0(t-r_j)
\end{equation}
we have
\begin{eqnarray*}
   \left(izI - \sum\limits_{j=0}^{N}A_{j,0}e^{-izr_j}\right)\hat{G}_0(z)&=&\Delta_{L_0}(iz)\hat{G}_0(z)\\
                                                                        &=&I,
\end{eqnarray*}
hence
\begin{equation}\label{g1}
G_0'(t)-\sum_{j=0}^NA_{j,0}G_0(t-r_j)=\delta(t)I
\end{equation}
with $\delta$  the Dirac delta  distribution.

 (\ref{g1}) implies that $G_0$ is an absolutely continuous function for all $t\neq0$, and it satisfies 
\begin{equation}
    G'_0(t) = \sum\limits_{j=1}^{N}A_{j,0}G_0(t+r_j) \ \ \mbox{ a. e.}\nonumber
\end{equation}
Also the function $G_0$ possess left- and right-hand limits $G_0(0-)$ and $G_0(0+)$, and there is a jump discontinuity 
\begin{equation}
    G_0(0+)-G_0(0-)=I.\nonumber
\end{equation}
\item[\textbf{Step 2}]\emph{  There are positive constant  $K, a_0>0$ such that $|G_0(t)|\leq Ke^{-a_0|t|}$. }
 \\
 Let
$$S(z):=-\sum\limits_{j=0}^NA_{j,0}e^{-r_jz},$$ and note that it is uniformly bounded in any vertical strip $$L=\{z\in \C/ \quad |\Re(z)|\leq a_0<1\},$$ for a fixed $0<a_0<1$.
This implies that, for $z\in L$ with $|\Im(z)|\gg 1$, the  matrix
\begin{equation}
\left|\di\frac{I-S(z)}{z+1}\right|<1,
\end{equation}                                
and in consequence we can write 

\begin{eqnarray*}
\Delta_{L_0}^{-1}(z)&=&(zI+S(z))^{-1}\\
                     &=&[(z+1)I+S(z)-I]^{-1}\\
                     &=&\di\frac{1}{z+1}\left[I-\di\frac{I-S(z)}{z+1}\right]^{-1}\\
                     &=&\di\frac{1}{z+1}\left\{I+\di\frac{I-S(z)}{z+1}+\di\frac{[I-S(z)]^2}{(z+1)^2}+\di\frac{[I-S(z)]^3}{(z+1)^3}+\cdots\right\}\\
                     &=&\di\frac{1}{z+1}I+\di\frac{I-S(z)}{(z+1)^2}+\di\frac{[I-S(z)]^2}{(z+1)^3}+\di\frac{[I-S(z)]^3}{(z+1)^4}+\cdots
\end{eqnarray*} 
This allows to write 
\begin{equation}\label{lambdainv}
\Delta_{L_0}^{-1}(z)=\di\frac{1}{z+1}I+R(z),
\end{equation}
where the  function $R(z):=\di\frac{I-S(z)}{(z+1)^2}+\di\frac{[I-S(z)]^2}{(z+1)^3}+\di\frac{[I-S(z)]^3}{(z+1)^4}+\cdots$ is analytic on $L$ and satisfies  $R(z)=O(|\Im(z)|^2)$ when $|\Im(z)|\to\infty$.  \\
Replacing (\ref{lambdainv}) into (\ref{G0}) 
\begin{eqnarray*}
G_0(t)&=&\frac{1}{2\pi}\int_{\mathbb{R}} e^{itz}\Delta^{-1}_{L_0}(iz)dz\\
      &=&\frac{1}{2\pi}\int_{\mathbb{R}} e^{itz} \left[\di\frac{1}{iz+1}I+R(iz)\right]dz\\
      &=&\frac{1}{2\pi}\int_{\mathbb{R}} e^{itz} \left[\di\frac{1}{iz+1}I\right]dz+\frac{1}{2\pi}\int_{\mathbb{R}} e^{itz} R(iz)dz\\
      &=&\frac{1}{2\pi}\int_{\mathbb{R}} e^{itz} \left[\di\frac{1}{iz+1}I\right]dz+\frac{1}{2\pi}\int_{\mathbb{R}} e^{itz} R(iz)dz  
\end{eqnarray*}

If $t>0$, the first last term is the inverse Fourier Transformation of   $$E(t):=\left\{\begin{array}{ll}0 &\mbox{if }t<0,\\ e^{-t}I &\mbox{if } t\geq 0. \end{array}\right.$$  In the second term we can change the path of integration, due to the analy\-ticity of $R(z)$ on $L$, from the line $\Re(z)=0$ to  the line $\Re(z)=-a_0$. It allows to write

\begin{equation}
G_0(t)=E(t)+\frac{e^{-a_0t}}{2\pi}\int_{\mathbb{R}} e^{itz} R(-a_0+iz)dz,
\end{equation}
showing that $|G_0(t)|\leq Ke^{-a_0t}$ for $t>0$, and some $K>0$. \\ For $t<0$ it is sufficiently change the path by $\Re(z)=a_0$ and use
 $$E_2(t)=\left\{\begin{array}{ll}0 &\mbox{if }t>0,\\ e^{t}I &\mbox{if } t\leq 0 \end{array}\right.$$ in spite of $E(t)$,  completing the proof of (\ref{G0-estimate}).  
\\
\item[\textbf{Step 3}]\emph{Given $h\in L^p_\mu$, the function $x(t):=\int_\R G_0(t-s)h(s)d\mu(s)$ is in $W^{1,p}_\mu$ and in consequence $\Lambda_{L_0}: W^{1,p}_\mu\to L^p_\mu$ is onto.}\\

We start checking $x\in L^p_{\mu}$. Let $\tilde{h}(z):=h(z)(1+z^2)^{-1/p}$

\begin{eqnarray}
    |x(t)| &\leq & \int_\R |G_0(t-z)||\tilde{h}(z)|(1+z^2)^{-1/q}dz\nonumber \\
    &\leq & \int_\R |G_0(t-z)||\tilde{h}(z)|dz  \nonumber \\
    &=&(|G_0|*|\tilde{h}|)(t).\nonumber   
\end{eqnarray}
Hence, using  Young's inequality we obtain
\begin{eqnarray}
\left(\int_\R|x(t)|^p(1+t^2)^{-1}dt\right)^{1/p} 
&\leq& \left(\int_\R((|G_0|*|\tilde{h}|)(t))^p(1+t^2)^{-1}dt\right)^{1/p}\nonumber\\
&\leq & \left(\int_\R((|G_0|*|\tilde{h}|)(t))^pdt\right)^{1/p}\nonumber\\
&\leq & \|G_0\|_{L^1}\|h\|_{L^p_\mu}\nonumber
\end{eqnarray}

impliying  $x\in L^p_{\mu}$.  In order to show that $x\in W^{1,p}_\mu$, let $\varphi\in C^{\infty}_c(\R)$ a test function and  we interprete $x$ as a regular generalized function as

\begin{equation}\label{rgeneralized}
T_x(\varphi)=\int_\R x(t)\varphi(t)d\mu,
\end{equation} 
wich implies that its generalized derivative is given by
\begin{equation}\label{dgeneralized}
T'_x(\varphi)=-\int_\R x(t)[\varphi(t)k(t)+\varphi'(t)]d\mu.
\end{equation}
In order to show that $x$ satisfies $\Lambda_{L_0}x=h$, we calculate

\begin{eqnarray}
&&\int_\R \varphi(t)\sum_{j=0}^NA_{j,0}x(t-r_j) d\mu(t)\nonumber\\
&=&\int_\R \varphi(t) \int_\R\left[\sum_{j=0}^NA_{j,0}G_0(t-r_j-s)\right]h (s)d\mu(s)\  d\mu(t)\nonumber\\
&=&\int_\R\varphi(t) \int_\R G_0'(t-s)h(s)d\mu(s)\  d\mu(t)\nonumber\\
&=&\int_\R h(s)w(s) \int_\R \varphi(t)[G_0'(t-s)w(t)]  dt\ ds\nonumber\\
&=&\int_\R h(s) w(s) \left[-\varphi(s)I- \int_\R \varphi(t)w'(t) G(t-s) dt+\int_\R \varphi(t) d_t(G(t-s)w(t))\right]ds\nonumber\\
&=&-\int_\R h(s)\varphi(s)d\mu(s)-\int_\R h(s)w(s) \int_\R \varphi(t)w'(t)G(t-s)dt  ds\nonumber\\
&&\quad+\int_\R h(s)w(s)\int_\R \varphi(t) d_t(G(t-s)w(t)) ds\nonumber\\
&=&-\int_\R h(s)\varphi(s)d\mu(s)-\int_\R \varphi(t) k(t) w(t)\left[\int_\R h(s) G(t-s)w(s)ds\right] dt\nonumber\\
&&\quad+\int_\R h(s)w(s)\left[-\int_\R\varphi'(t)G(t-s)w(t)dt\right] ds\nonumber\\
&=&-\int_\R h(s)\varphi(s)d\mu(s)-\int_\R \varphi(t) k(t) x(t) d\mu(t)\nonumber \\ 
&&\quad-\int_\R \varphi'(t)\left[\int_\R G(t-s)h(s)w(s)ds\right] w(t)dt\nonumber\\
&=&-\int_\R h(s)\varphi(s)d\mu(s)-\int_\R [\varphi(t) k(t) +\varphi'(t)]x(t)d\mu(t) \nonumber
\end{eqnarray}
it implies that $x$ satisfies the generalized version of $\Lambda_{L_0}x=h$

\begin{equation}\label{weekeq}
\int_\R \varphi(t)\left[\sum_{j=0}^NA_{j,0}x(t-r_j)+h(t)\right] d\mu(t)=-\int_\R [\varphi(t) k(t) +\varphi'(t)]x(t)d\mu(t).
\end{equation}

Let   $g(t)=\sum_{j=0}^{N}A_{j,0}x(t-r_j)+h(t)\in L^p_{\mu}$. If we define $u(t)= \int_{t_0}^tg(s)ds$ for some $t_0\in\R$,  then $u$ is absolutely continuous and integrating by parts and using (\ref{weekeq})
\begin{eqnarray}
T'_u(\varphi)&=&-\int_\R u(t)[\varphi'(t)+k(t)\varphi(t)]d\mu\nonumber\\
             &=&\int_\R u'(t)\varphi(t)d\mu \nonumber\\
             &=&T'_x(\varphi).\nonumber
\end{eqnarray}
Hence 
\begin{eqnarray}
0&=&T'_u(\varphi)-T'_x(\varphi)\nonumber\\
 &=&\int_\R[x(t)-u(t)][\varphi'(t)+k(t)\varphi(t)]d\mu\nonumber\\
 &=&\int_\R[x(t)-u(t)][\varphi(t) w(t)]'dt\nonumber\\
 &=&\int_\R[x(t)-u(t)]\phi'(t) dt,\nonumber
\end{eqnarray}
for all $\phi\in C^{\infty}_c(\R)$. We conclude that $x(t)=u(t)+C$, for some $C\in\R$ and in consequence $x$ is absolutely continuous function and satisfying, a.e. the equation

\begin{equation}\label{weekeq2}
x'(t)=\sum_{j=0}^NA_{j,0}x(t-r_j)+h(t).
\end{equation}

 This implies that $x\in W^{1,p}_{\mu}$ and that $\Lambda_{L_0}x=h$ and so $\Lambda_{L_0}:W_{\mu}^{1,p}\to L_{\mu}^p$ is onto. 
 
 \item[\textbf{Step 4}] {The operator $\Lambda_{L_0}:W_{\mu}^{1,p}\to L_{\mu}^p$ is one-to-one.}\\
Indeed, we suppose that $\Lambda_{L_0}x=0$ for some $x\in W_{\mu}^{1,p}$, i.e $x$ satisfies in $\R$
\begin{equation}\label{homogenea}
x'(t)=\sum_{j=0}^NA_{j,0}x(t-r_j).
\end{equation}
 
   Interpreting $x$ as tempered generalized function and appliying the  Fourier transform on (\ref{homogenea}), we obtain that    $\hat{x}$ satisfies the equality $\Delta_{L_0}(iz)\hat{x}(z)=0$. Now from (\ref{hip-cond}) we concluded that $\hat{x}$ is the zero generalized function and hence $x$ is the zero function in $W^{1,p}$ completing the proof of the one-to-one property. 
\end{enumerate}
Finally, from the steps 3 and 4 is clear that $\displaystyle (\Lambda_0^{-1}h)(t)=\int_\R G(t-s)h(s)d\mu$ concluding the proof.
\end{proof}
\newpage

\begin{lemma}\label{Prop-G(t,z)}
Let
\begin{equation}
x'(t)=L(t)x_t=L_0x_t+M(t)x_t
\end{equation}
with $x'(t)=L_0x_t$ an  hyperbolic system. There are positive constants  $\epsilon, K$ and $a$ such that 
\begin{enumerate}
\item[(i)] If \begin{equation}\label{smallM}
\|M(t)\|\leq \varepsilon,\quad \mbox{for all }t\in\R,
\end{equation} 
then the asociated operator $\Lambda_{L}: W_{\mu}^{1,p}\to L_{\mu}^p$ is an isomorphism for $1\leq p\leq \infty$. 
\item[(ii)] There exist a function $G:\mathbb{R}^2\to \mathbb{C}^d$  satisfying 
\begin{equation}\label{Gestimation}
    |G(t,z)|\leq Ke^{-a|t-z|}, \ \  \mbox{ for all }(t,z)\in\R^2,
\end{equation}
and 
\begin{equation}\label{Lambda_Linversa}
    (\Lambda_{L}^{-1}h)(t)=\int_\R G(t,z)h(z)d\mu(z).
\end{equation}
\end{enumerate}

\end{lemma}

\begin{proof}
Writing  $\Lambda_L=\Lambda_{L_0}-M$, and using that $\Lambda_{L_0}$ has inverse by Theorem (\ref{T-Green-L_0}) we choose $\varepsilon<\|\Lambda^{-1}_{L_0}\|$, wich implies $\|M\Lambda_{L_0}^{-1}\|_{\mathcal{L}(L_{\mu}^p, L_{\mu}^p)}<1$, and in consequence $ \Lambda^{-1}_{L}:L_{\mu}^p \to W_{\mu}^{1,p}$ exists due to the convergence of its  Neumanns serie

\begin{eqnarray} 
\Lambda^{-1}_{L}&=&\Lambda^{-1}_{L_0}(I-M\Lambda^{-1}_{L_0})^{-1}\nonumber\\
         &=& \Lambda^{-1}_{L_0}\sum\limits_{j=0}^{\infty}(M\Lambda^{-1}_{L_0})^{j},\label{neumann}
\end{eqnarray}
proving $(i)$.\\

In order to prove (\ref{Gestimation}), let $G_0$ the Green's function associated to $\Lambda_{L_0}$ in the Theorem (\ref{T-Green-L_0}) and let $h\in L^p_\mu$, hence
\begin{eqnarray}
    (M\Lambda_{L_0}^{-1}h)(t) &=& M\left(\int_{\R}G_0(t-z) h(z)\ d\mu(z)\right)\nonumber\\
                              &=&\sum\limits_{j=0}^NC_j(t)\int_{\R}G_0(t-r_j-z) h(z)\ d\mu(z)\nonumber\\
                              &=&\int_{\R}\left[\sum\limits_{j=0}^NC_j(t)G_0(t-r_j-z)\right] h(z)\ d\mu(z)\nonumber\\
                              &=&\int_{\R}\Gamma(t,z) h(z)\ d\mu(z)\label{gamma},  
\end{eqnarray}
Using (\ref{G0-estimate}) with constants $a_0$, $K_0$, we have that the function
$$\Gamma(t,z)=\sum\limits_{j=0}^{N}C_j(t)G_0(t-r_j-z)$$ satisfies the following estimation

\begin{equation}\label{Gamma-estimate}
    |\Gamma(t,z)|\leq \sum\limits_{j=0}^N |C_j(t)|K_0 e^{a_0|r_j|}e^{-a_0|t-z|}\leq K_{1}e^{-a_0|t-z|},
\end{equation}
where 
$$K_{1}=K_0\varepsilon e^{a_0r_N}.$$

If we write
$$((M\Lambda_{L_0}^{-1})^{j}h)(t)=\int_\R\Gamma_j(t,z)h(z)dz,$$
with 
\begin{equation}
    \Gamma_j(t,z):=\left\{\begin{array}{ll} \Gamma(t,z), & \mbox{if }j=1,\\ & \\
                                         \displaystyle\int_\R\Gamma(t,s)\Gamma_{j-1}(s,z)ds, &\mbox{if }\  j\geq 2,
                           \end{array}\right.\nonumber              
\end{equation}
from (\ref{neumann}) our candidate to be a Green's function satisfying $(ii)$ is

\begin{equation}\label{green-L}
    G(t,z) = G_0(t-z) + \int_{\R}G_0(t-s)\left(\sum\limits_{j=1}^{\infty}\Gamma_j(s,z)\right)ds.
\end{equation}

Now we want to  show the exponential estimation for $G(t,z)$. For that, let $\psi(t):=K_1e^{-a_0|t|}$ and 
$$\psi^{*j}(t):=\left\{\begin{array}{ll}\psi(t),&\mbox{ if } j=1\\ (\psi * \psi^{*(j-1)})(t),&\mbox{ if }j\geq 2.\end{array}\right.
$$
 Appliying  (\ref{Gamma-estimate}) inductively we have 
\begin{equation}
|\Gamma_j(t,z)|\leq \psi^{*j}(t-z). 
\end{equation}
Assuming $\varepsilon<\di\frac{a_0e^{-a_0r_N}}{2 K_0}$, which implies $K_1<\di\frac{a_0}{2}$,  we have  by  Lemma 5.1 in \cite{Mallet-Paret}, 
\begin{equation} 
\sum\limits_{j=1}^{\infty}|\Gamma_j(t,z)|\leq K_2e^{-a_1|t-z|}
\end{equation}

where  $a_1=\sqrt{a_0^2-2a_0K_1}$ and $K_2=\di\frac{K_1}{a_1}$ . \\
From (\ref{green-L}) and due to $a_1=a_0\sqrt{1-2K_1/a_0}<a_0$, there exist $K_3$ such that  
\begin{eqnarray}
|G(t,z)|&\leq &K_0e^{-a_0|t-z|}+\int_\R K_0K_2e^{-a_0|t-s|-a_1|s-z|}ds\nonumber\\
        &\leq& K_3e^{-a_1|t-z|},   \nonumber
\end{eqnarray}
for all $(t,z)\in \R^2$, showing (\ref{Gestimation}).
\\ 
Finally, to show the representation (\ref{Lambda_Linversa}), let any $h\in L^p_{\mu}$, and the corresponding $x\in W^{1,p}_\mu$ such that   $\Lambda_L x=h$. We know that the sequence $(x_k)_{k\in\N}$ defined by
\begin{equation}
    x_k:=\Lambda^{-1}_{L_0}\sum\limits_{j=0}^{k}(M\Lambda^{-1}_{L_0})^j h\nonumber
\end{equation}
meets that  $x_k\to x$ in $W_{\mu}^{1,p}$ and is easy to see that  
\begin{equation}
    x_k(t) = \int_{\R} G_k(t,z)h(z)d\mu(z),\nonumber
\end{equation}

where
\begin{equation}
    G_k(t,z) = G_0(t-z) + \int_\R G_0(t-s)\Big(\sum\limits_{j=1}^{k} \Gamma_j(s,z)\Big)ds.\nonumber
\end{equation}
Hence, we want to show that $x(t)=y(t)$, where 
\begin{equation}
y(t)=\int_\R G(t,z)h(z) d\mu(z).
\end{equation} 
In order to do that,  note that

\begin{eqnarray}
    |x_k(t)-y(t)|&=&\left|\int_{\R}\left[G_k(t,z)-G(t,z)\right] h(z)d\mu(z)\right|\nonumber\\
                 &\leq &\int_{\R}\int_{\R} |G_0(t-s)|\sum\limits_{j=k+1}^\infty\psi^{*j}(s-z)\ ds |h(z)|d\mu(z)\nonumber\\
                 &= &\int_{\R}\int_{\R} |G_0(t-z-s)|\sum\limits_{j=k+1}^\infty\psi^{*j}(s)\ ds |h(z)|d\mu(z)\nonumber\\
                 &= &\int_{\R}R_k(t-z) |h(z)|d\mu(z)\nonumber\\
                 &=& \int_{\R}R_k(t-z) |h(z)|(1+z^2)^{-1/p}(1+z^2)^{-1/q}dz\nonumber\\
                  &\leq & \int_{\R}R_k(t-z) |h(z)|(1+z^2)^{-1/p}dz\nonumber\\
                  &= & \int_{\R}R_k(t-z) |\tilde{h}(z)|dz\nonumber\\
                  &=&  (R_k*|\tilde{h}|)(t),\nonumber
\end{eqnarray}
where $R_k(t)=\left(G_0*\sum\limits_{j=k+1}^\infty\psi^{*j}\right)(t)$.  Hence
\begin{eqnarray}
\int_{\R}|x_k(t)-y(t)|^p\ d\mu(t)&\leq & \int_\R [(R_k*|\tilde{h}|)(t)]^p(1+t^2)^{-1}dt\nonumber\\
                                 &\leq & \int_\R [(R_k*|\tilde{h}|)(t)|]^p dt\nonumber\\
                                 &\leq & \|R_k\|_{L^1}^p\| \tilde{h}\|_{L^p}^p\nonumber\\
                                 &=&\|R_k\|_{L^1}^p\| h\|_{L_\mu ^p}^p\nonumber
\end{eqnarray}
By the last estimation we only need to show that $R_k\to 0$ in $L_1$ to conclude $x_k\to y$ in $L_\mu^p$.  For that, applying one more time  Young's inequality is easy to see that

\begin{equation}
\|R_k\|_{L^1}\leq \|G_0\|_{L^1}\sum_{j=k+1}^\infty\|\psi\|_{L^1}^j\nonumber
\end{equation}
showing that $R_k\to 0$ when $k\to \infty$ and in consequence  $x_k\to y$ in $L_{\mu}^p$ and thus $x=y$.
\end{proof}

\begin{lemma}\label{P-cotas-x}
Assuming that the equation \begin{equation}\label{eq3}
x'(t)=L(t)x_t
\end{equation}
 is asymptotically hyperbolic at $+\infty$ (respectively at $-\infty$). Then 
 \begin{enumerate}
 \item There exists  constants $a$, $K_1$,  and $K_2$, all positives, such that for  $x\in W_{\mu}^{1,p}$, and $h\in L_{\mu}^p$, satisfying  $\Lambda_{L}x=h$   we have 
\begin{eqnarray}
    |x(t)| &\leq& K_1e^{-a|t|}\|x\|_{L^{\infty}_{\mu}} + K_1\int_{\R}e^{-a|t-z|}|h(z)|d\mu(z)\label{x-cota1}\\
    &\leq & K_1e^{-a|t|}\|x\|_{L^{\infty}_{\mu}} + K_2\|h\|_{L_{\mu}^p}\label{x-cota2}
\end{eqnarray}
for $t\geq 0$ (respectively $t\leq 0$). 
\item If (\ref{eq3}) is asymptotically hyperbolic in both $\pm\infty$, then (\ref{x-cota1}) and (\ref{x-cota2}) are valid on all $\R$ and there exists $K_3>0$ such that
\begin{equation}\label{x-cota3}
    \|x\|_{W_{\mu}^{1,p}} \leq K_3(\|x\|_{L_{\mu}^{\infty}}+\|h\|_{L^p_{\mu}}).
\end{equation}
 \end{enumerate}
\end{lemma}

\begin{proof} We assume that (\ref{eq3})  is asymptotically hyperbolic at $+\infty$, the case at $-\infty$ is analogous, and also from both cases we deduce (\ref{x-cota1}) and (\ref{x-cota2}) on all $\R$. Hence we assume $t\geq 0$.\\
  We write $L(t)=L_0+M(t)$, where $\|M(t)\|\to 0$ as $t\to +\infty$.
   Let $\varepsilon>0$ as in Lemma (\ref{Prop-G(t,z)}) and  let  $\tau>0$ such that $\|M(t)\|\leq \varepsilon$ for all $t\geq \tau$. Now we define 
 
\begin{equation}\label{alpha}
    \alpha(t)=\left\{\begin{array}{lr}
        0, & t<\tau \\
        1, & t\geq\tau
    \end{array}\right.,
\end{equation}

\begin{equation}
    L_1(t)=L_0+\alpha(t)M(t),\nonumber 
    \end{equation}

and \begin{equation}
    M_1(t)=(1-\alpha(t))M(t).\nonumber
\end{equation}
Note that   $M_1(t)=0$ if $t\geq \tau$.
With these definitions, for $x$ and $h$ satisfying $\Lambda_Lx=h$, we rewrite (\ref{eq3})
\begin{equation}
    x'(t) = L_{1}(t)x_t + M_1(t)x_t + h(t).\nonumber
\end{equation}
Applying Lemma (\ref{Prop-G(t,z)}) to $x'(t) = L_1(t)x_t$, and note that   $\|\alpha(t)M(t)\|\leq \varepsilon$ for all $t\in \mathbb{R}$, and for the same $\varepsilon$, because it only depends of $L_0$, there exist  $G_1:\mathbb{R}^2\to \mathbb{C}^d$ a Green's function for $\Lambda_{L_1}$ and constants $K,a$ such that

\begin{eqnarray}
x(t) &=& \int_{\R} G_1(t,z)(M_1(z)x_z + h(z))d\mu(z)\nonumber\\
&=& \int_{-\infty}^{\tau} G_1(t,z)M_1(z)x_z(1+z^2)^{-1}dz + \int_{\R} G_1(t,z)h(z)d\mu(z)\nonumber
\end{eqnarray}

Hence, for any $t\in \R$, we have that

\begin{eqnarray}
    |x(t)| &\leq & K\int_{-\infty}^{\tau} e^{-a|t-z|}|M_1(z)x_z|(1+z^2)^{-1}dz + K\int_{\R} e^{-a|t-z|}|h(z)|d\mu(z)\nonumber\nonumber\\
           &\leq & K\int_{-\infty}^{\tau} e^{-a|t-z|}\|M_1(z)\| \|x\|_{L_\mu^{\infty}} dz + K\int_{\R} e^{-a|t-z|}|h(z)|d\mu(z)\nonumber\nonumber\\
           &\leq & K\left\{\sup\limits_{z\in\R}\|M_1(z)\|\right\} \|x\|_{L_\mu^{\infty}}\int_{-\infty}^{\tau} e^{-a|t-z|} dz + K\int_{\R} e^{-a|t-z|}|h(z)|d\mu(z).\nonumber
\end{eqnarray}

Is easy to see that, for $t\geq \tau$,  $$ \int_{-\infty}^{\tau} e^{-a|t-z|} dz\leq Ce^{-a|t|},$$ with $C>0$ and in consequence for some $C_2\geq C$ this inequality  remais for $t\geq 0$. Hence, from 
\begin{equation}
|x(t)|\leq K\left\{\sup\limits_{z\in\R}\|M_1(z)\|\right\}C_2 \|x\|_{L_\mu^{\infty}} e^{-a|t|} + K\int_{\R} e^{-a|t-z|}|h(z)|d\mu(z)
\end{equation}
 we can  choose $K_1>0$ in a way that the estimate in (\ref{x-cota1}) holds for all $t\geq0$. The choice of $K_2$ and the proof of inequality (\ref{x-cota2}) follow from
 
 \begin{eqnarray}
 K\int_{\R} e^{-a|t-z|}|h(z)|d\mu(z)&\leq & K\int_{\R} e^{-a|t-z|}|h(z)|(1+z^2)^{-1/p}dz\nonumber\\
                                    &\leq & K\left[\int_{\R} e^{-aq|t-z|}dz\right]^{1/q}\left[\int_\R|h(z)|^p(1+z^2)^{-1}dz\right]^{1/p}\nonumber\\
                                    &=& K_2\|h\|_{L^p_\mu},\nonumber
 \end{eqnarray}
 with $K_2=K\left[\int_{\R} e^{-aq|z|}dz\right]^{1/q}$.

In order to obtain (\ref{x-cota3}), we assume valid (\ref{x-cota1}) for all $t\in\R$, and from it, calling $\gamma(t):=e^{-a|t|}$, $\tilde{h}(t):=h(t)(1+t^2)^{-1/p}$ 

\begin{eqnarray}
|x(t)|(1+t^2)^{-1/p}&\leq &K_1\|x\|_{L_\mu^\infty}\gamma (t)(1+t^2)^{-1/p} +K_1(\gamma*|\tilde{h}|)(t)(1+t^2)^{-1/q}\nonumber\\
                     &\leq &K_1\|x\|_{L_\mu^\infty}\gamma (t)(1+t^2)^{-1/p} +K_1(\gamma*|\tilde{h}|)(t),\nonumber
\end{eqnarray}
hence, by Minkowski and Young's inequality

\begin{eqnarray}
\left(\int_\R|x(t)|^pd\mu(t)\right)^{1/p}&\leq &K_1\|x\|_{L_\mu^\infty}\left(\int_\R\gamma^p(t)(1+t^2)^{-1} dt\right)^{1/p}+K_1\left(\int_\R(\gamma*|\tilde{h}|)^p(t)dt\right)^{1/p}\nonumber\\
&\leq & K_1\|\gamma\|_{L^p_\mu}\|x\|_{L^\infty_\mu}+K_1\|\gamma\|_{L^1}\|h\|_{L_\mu^p},\nonumber
\end{eqnarray}
and letting $k_1=K_1\max\{\|\gamma\|_{L^p_\mu},\|\gamma\|_{L^1}\}$ we arrive to

\begin{equation}\label{x-cota4}
 \|x\|_{L^p_\mu} \leq k_1\left(\|x\|_{L^{\infty}_{\mu}}+ \|h\|_{L_\mu^p}\right).
\end{equation}

 Now, from the differential equation (\ref{eq2}) and the boundedness of the coefficients $A_j$ for $t\in\mathbb{R}$ we concluede  that $x'\in L^p_{\mu}$, and
\begin{equation}\label{x-cota5}
    \|x'\|_{L^p_{\mu}} \leq k_2(\|x\|_{L^p_{\mu}} + \|h\|_{L^p_{\mu}}),
\end{equation}
for some $k_2$ depending only on $L$. \\
Finally, we  obtain (\ref{x-cota3}) add (\ref{x-cota4}) and (\ref{x-cota5}) and choose $K_3=\max\{k_1,k_2\}$.
\end{proof} 
The following properties have proofs very close to the  case without weight  showed in \cite{Mallet-Paret}, we summary these in the following lemma and  omit the proofs.
\begin{lemma}\label{summary} Assuming that the equation (\ref{eq0}) is asymptotically hyperbolic and $1\leq p\leq \infty$.
\begin{enumerate}
\item  If there are  sequences $x_n\in W_{\mu}^{1,p}$ and $h_n\in L^p_{\mu}$ with $\Lambda_L x_n=h_n$, such that the sequence $x_n$  is bounded in $W_{\mu}^{1,p}$ and  $h_n\to h_{*}$ in $L^p_{\mu}$. Then there exists a subsequence $x_{n_k}\to x_{*}$ converging in $W_{\mu}^{1,p}$ to some element $x_{*}$, with necessarily $\Lambda_L x_{*}=h_{*}$.
\item 
The kernel $\mathcal{K}^p_L$  of the operator $\Lambda_L$ is independent of $p$ and finite-dimensional.
\item 
For $1\leq p\leq \infty$ the range $\mathcal{R}_L^p\subseteq L_{\mu}^p$ of $\Lambda_L$ is closed.
\end{enumerate}
\end{lemma}


 Finally, we prove  Theorem \ref{T-Fredholm}

\begin{proof} As we are assuming that the equation (\ref{eq0}) is asymptotically hyperbolic, from the last lemma  we have that the kernel $\mathcal{K}_L:=\mathcal{K}_L^p$ of $\Lambda_L$ is independent of $p$,   and that $\dim \mathcal{K}_L<\infty$. Moreover,
by $(P_1)$, $(P_2)$ and (\ref{B_j}), the adjoint equation (\ref{adjunta_of_L}) is also asymptotically hyperbolic and so applying the same lemma we have that the subspace $\mathcal{K}_{L^{*}}\subseteq W^{1,q}_{\mu}$, with $\frac{1}{p}+\frac{1}{q}=1$, is finite-dimensional and independent of $q$, showing (\ref{T-Fredholm})-(1) and (\ref{T-Fredholm})-(2).
\\
In order to show (\ref{T-Fredholm})-(3), we consider the set $\mathcal{S}_L^p\subseteq L_{\mu}^p$ as the subspace 
\begin{equation}
    \mathcal{S}_L^p = \Big\{h\in L_{\mu}^p : \int_\R\overline{y(t)}h(t)d\mu=0, \mbox{ for all } y\in\mathcal{K}_{L^{*}} \Big\},\nonumber
\end{equation}
 hence $\mbox{Codim }\mathcal{S}_L^p = \dim \mathcal{K}_{L^{*}}<\infty$ in $L^p_{\mu}$. We want to show $\mathcal{R}_L^p=\mathcal{S}_L^p$. To this end, let any $h=\Lambda_L x\in \mathcal{R}_L^p$, with $x\in W^{1,p}_{\mu}$, and $y\in \mathcal{K}_{L^{*}}$.  From (\ref{igualdad-adj}) we have that
\begin{equation}
    \int_\R \overline{y(t)}h(t)d\mu = \int_\R\overline{y(t)}(\Lambda_L x)(t)d\mu = \int_\R\overline{(\Lambda^{*}_L y)(t)}x(t)d\mu = 0,\nonumber
\end{equation}
 which implies that $h\in \mathcal{S}_L^p$ and so, $\mathcal{R}_L^p\subseteq\mathcal{S}_L^p$. Conversely, assume first that $p\neq \infty$ and  $\mathcal{R}_L^p\neq\mathcal{S}_L^p$. In these conditions, and as $\mathcal{R}_L^p$ is closed by (\ref{summary})-(3),   there exists a nontrivial linear continuous functional on the space $L^p_{\mu}$ vanishing identically on $\mathcal{R}_L^p$, but not identically on $\mathcal{S}_L^p$, also it is represented by some $y\in L^q_{\mu}$, that is
\begin{equation}\label{cont1}
    \int_\R\overline{y(t)}(\Lambda_L x)(t)d\mu = 0, \qquad \mbox{ for all } x\in W^{1,p}_{\mu} 
\end{equation}
and
\begin{equation}\label{cont2}
    \int_\R\overline{y(t)}h(t)d\mu \neq 0 \qquad \mbox{ for some } h\in\mathcal{S}_L^p.
\end{equation}

Let $\xi:\mathbb{R}\to \mathbb{C}^d$ a $C^{\infty}$ function  with compact support, and set $x(t)=\overline{\xi(t)}$. Then, from (\ref{cont1}) we have 
\begin{eqnarray}
    0 &=& \int_\R y(t)\overline{(\Lambda_L \overline{\xi})(t)}d\mu\nonumber\\
    &=& \int_\R y(t)\xi'(t)d\mu - \sum\limits_{j=0}^{N}\int_\R y(t)\overline{A_j(t)}\xi(t-r_j)d\mu\nonumber\\
    &=& \int_\R y(t)\xi'(t)d\mu + \int_\R y(t)\xi(t)k(t)d\mu\nonumber\\
    && -\int_\R y(t)\overline{[k(t)I+A_0(t)]}\xi(t)d\mu -\sum\limits_{j=1}^{N}\int_\R y(t+r_j)\overline{M_j^+(t)A_j(t+r_j)}\xi(t)d\mu\nonumber\\
    &=& \int_\R y(t)[\xi'(t) +k(t)\xi(t)]d\mu - \int_\R \left[\sum\limits_{j=0}^{N}B^{*}_j(t) y(t+r_j) \right]\xi(t)d\mu \nonumber
\end{eqnarray}

Following the same argument as in the proof of Theorem (\ref{T-Green-L_0})-Step 3, the above formula  says that the generalized derivative of  $y$, 
defined as in (\ref{dgeneralized}), satisfies the adjoin equation (\ref{adjunta_of_L}), and in consequence   $y\in W^{1,q}_{\mu}$ and  $y\in \mathcal{K}_{L^{*}}$. This prove that $\mathcal{R}_L^p=\mathcal{S}_L^p$.\\
The case $p=\infty$ is proved in \cite{Mallet-Paret} Theorem A, and it is analogous for our  spaces.  
\end{proof}

\section{Proof of Theorem \ref{T-E-D}}
We start noting  that from Theorem \ref{T-Fredholm} and Lemma (\ref{summary}),  the operator $\Lambda_L: W^{1,p}_{\mu}\to L^p_{\mu}$ is a Fredholm operator, closed and  with close range. If we assume that $\mbox{Ind}(\Lambda_L)=\dim(\mathcal{K}_L)=0$, there exists its inverse $\Lambda_L^{-1}: L^p_{\mu}\to W^{1,p}_{\mu}$ and it must be bounded. Onward, we call
\begin{eqnarray}
 r&:=&r_{N},\ \mbox{ the longest delay},\nonumber\\
 \beta&:=&\sup\limits_{t\in\R}\|L(t)\|,\nonumber\\
 \gamma_0&:=&\|\Lambda^{-1}_L\|\beta+1.\nonumber
\end{eqnarray}
Also, given any $s\in\R$ and $\phi\in C$ we define two functions:
\begin{enumerate}
\item $\psi:[s-r,\infty)\to \R^n$, 
\begin{equation}\label{psi-function}
\psi(t)=\left\{\begin{array}{ll}\phi(t-s),&\mbox{ if } t\in [s-r,s],\\
                                                                   \phi(0),&\mbox{ if } t\in (s,+\infty)
                                                  \end{array}\right.
                                                  \end{equation}
\item $g:\R\to\R^n$, 
\begin{equation}\label{g-function}
g(t)=\left\{\begin{array}{ll}0,&\mbox{ if } t\in (-\infty, s),\\
                       L(t)\psi_t,&\mbox{ if } t\in [s,+\infty)
                                                  \end{array}\right.
                                                  \end{equation}

\end{enumerate}
Note that $g$ is measurable and  
\begin{equation}\label{g-sup}
\sup\limits_{t\in\mathbb{R}}|g(t)|\leq \beta\|\phi\|<+\infty,
\end{equation}
hence $g\in L_\mu^p$, with
 \begin{equation}\label{g-cotaLp}
 \|g\|_{L^p_\mu}\leq \sqrt[p]{\pi}\beta\|\phi\|
 \end{equation}

  In order to show the existence of an exponential dichotomy on all $\R$, we divide the proof in the lemmas below

\begin{lemma}\label{lemma-split}
For each $s\in \mathbb{R}$  let
\begin{equation}\label{E(s)}
    E(s):=\{\phi\in C : \sup\limits_{t\geq s}\|T(t,s)\phi\|<+\infty\}
\end{equation}
and 
\begin{equation}\label{F(s)}
\begin{array}{rcl}
    F(s) &:= &\{\phi\in C : \mbox{ there exists } v:(-\infty, s]\to \mathbb{R}^n \mbox{ continuous with }\\
    & & v_s=\phi, \quad \sup\limits_{t\leq s}|v(t)|<+\infty  \mbox{ and } v_t=T(t,\tau)v_{\tau} \mbox{ for } s\geq t\geq \tau\}.
    \end{array}
\end{equation}
Then, we have the split
\begin{equation}\label{split}
    C=E(s)\oplus F(s).
\end{equation}
\end{lemma}
\begin{proof}
Taking $s\in\mathbb{R}$, $\phi\in C$, and their associated functions  $\psi:[s-r, +\infty)\to \mathbb{R}^n$ and 
and $g:\mathbb{R}\to \mathbb{R}^n$, we have that there is some $v\in W_\mu^{1,p}$ such that $\Lambda_Lv=g$. Hence the function $u=v+\psi$ is a solution of the equation $x'(t)=L(t)x_t$, with $u_s=v_s+\phi$. That is,
\begin{equation}
    u_t=T(t,s)(v_s+\phi) \quad \mbox{ for } \quad t\geq s.\nonumber
\end{equation}

From Lemma  (\ref{P-cotas-x}), inequality (\ref{x-cota2}),  we obtain that there exists positive cons\-tants $K_1,K_2$ and $a$ such that 
\begin{eqnarray}
    |v(t)| &\leq & K_1e^{-a|t|}\|v\|_{L^{\infty}_{\mu}} + K_2\|g\|_{L^p_\mu}\nonumber
    \end{eqnarray}
   which implies that $\sup\limits_{t\in\mathbb{R}}|v(t)|<+\infty$ and so 
\begin{equation}
    \sup\limits_{t\geq s}|u(t)|\leq \sup\limits_{t\geq s}|v(t)| + \sup\limits_{t\geq s} |\psi(t)| <+\infty.\nonumber
\end{equation}
In conclusion,  $v_s+\phi\in E(s)$.\\
 By  another hand, since $g(t)=0$ for $t<s$, we have that $\Lambda_L v=0$ for $t<s$, that is, $v_t=T(t,\tau)v_\tau$, for $\tau\leq t\leq s$. Hence, $v_s\in F(s)$ and 
\begin{equation}
    \phi = (v_s + \phi) - v_s \in E(s) + F(s).\nonumber
\end{equation}
Now, let $\phi\in E(s)\cap F(s)$. Then, there exists a continuous function $v:(-\infty, s]\to \mathbb{R}^n$ with $v_s=\phi$ such that $\sup\limits_{t\leq s}|v(t)|<+\infty$ and $v_t=T(t,\tau)v_{\tau}$, for $\tau\leq t\leq s$. We define a function $u: \mathbb{R}\to \mathbb{R}^n$ by 
\begin{equation}
    u_t = \left\{ \begin{array}{ll}
        T(t,s)\phi, & \mbox{ if } t\geq s,  \\
        v_t, & \mbox{ if } t< s. 
    \end{array} \right.\nonumber
\end{equation}
Observe that $u$ is continuous and since $\phi\in E(s)$ we have that $\sup\limits_{t\in\mathbb{R}}|u(t)|<+\infty$, hence  $u\in L^p_{\mu}$ and it is easy to verify that $u_t=T(t,\tau)u_{\tau}$ for any  $\tau\leq t$, which implies that $\Lambda_L u=0$ and in consequence  $u\in W^{1,p}_{\mu}$. It follows from the invertibility of  $\Lambda_L$  that $u=0$ and so $\phi=v_s=u_s=0$. Therefore $E(s)\cap F(s)=\{0\}$.
\end{proof}

\begin{lemma}\label{lemma-T|F}
 Let
\begin{equation}\label{Projections}
    P(s), Q(s): C\to C
\end{equation}
be the projections associated with the splitting $C=E(s)\oplus F(s)$. Then, for each $t\geq s$, the linear operator 
\begin{equation}\label{T|F}
    T(t,s)|F(s) : F(s)\to F(t) 
\end{equation}
is onto and invertible.
\end{lemma}
\begin{proof}
Consider that $T(t,s)\phi=0$ for some $\phi\in F(s)$. Then there exists a continuous function $v:(-\infty, s]\to \mathbb{R}^n$ with $v_s=\phi$ such that $\sup\limits_{\tau\leq s}|v(\tau)|<+\infty$ and $v_{\tau}=T(\tau, \overline{\tau})v_{\overline{\tau}}$,  for $\overline{\tau}\leq \tau\leq s$. Defining a function $u:\mathbb{R}\to \mathbb{R}^n$ by
\begin{equation}
    u_{\tau} = \left\{ \begin{array}{ll}
    0, & \mbox{ if } t\leq \tau,  \\
        T(\tau,s)\phi, & \mbox{ if } s\leq \tau <t,  \\
        v_{\tau}, & \mbox{ if } \tau\leq s. 
    \end{array} \right.\nonumber
\end{equation}
hence $u$ is bounded and continuous on $\R$,  and in consequence $u\in L^p_\mu$. Moreover, we have that $u_{\tau}=T(\tau, \overline{\tau})u_{\overline{\tau}}$ for any real  $\tau\geq\overline{\tau}$. Thus, $\Lambda_L u=0$ and then $u\in W^{1,p}_{\mu}$. Again, the invertibility of $\Lambda_L$ implies that $u=0$ and hence $\phi=v_s=u_s=0$. This shows that the operator in (\ref{T|F}) is one-to-one.

Now, take $\phi\in F(t)$. Then, there exists a continuous function $v:(-\infty, t]\to \mathbb{R}^n$ with $v_t=\phi$ such that $\sup\limits_{\tau\leq t}|v(\tau)|<+\infty$ and $v_{\tau}=T(\tau,\overline{\tau})v_{\overline{\tau}}$ for $\overline{\tau}\leq \tau\leq t$. In particular,
\begin{equation}
    \phi=v_t=T(t,s)v_s\nonumber
\end{equation}
and since $v_s\in F(s)$, we obtain that the operator in (\ref{T|F}) is onto.
\end{proof}
In \cite{Barreira-valls-hyperbolicity}, Lemma 4.6, was proved that 
\begin{equation}
    P(t)T(t,s) = T(t,s)P(s) \quad \mbox{ for } t\geq s.\nonumber
\end{equation}
and the proof is still valid for our case, hence we assume it onward.
\begin{lemma}
There projection $P(s)$ satisfies  
\begin{equation}\label{P-cota}
    \|P(s)\|\leq \gamma_0,  \mbox{ for }  s\in\mathbb{R}.
\end{equation}
\end{lemma}

\begin{proof}
Taking $s\in\mathbb{R}$, $\phi\in C$, and  $v=\Lambda^{-1}_Lg$, we have 

\begin{eqnarray}
    \|P(s)\phi\| = \|v_s+\phi\| &\leq& \sup\limits_{t\in\mathbb{R}}|v (t)| + \|\phi\|\\
    &=& \sup\limits_{t\in\mathbb{R}}|(\Lambda^{-1}_L g)(t)| + \|\phi\|\\
    &\leq & \|\Lambda^{-1}_L\|\sup\limits_{t\in\mathbb{R}}|g(t)| + \|\phi\|\\
    &\leq &  (\|\Lambda^{-1}_L\|\beta + 1)\|\phi\|\\
    &=&\gamma_0\|\phi\|.
\end{eqnarray}
 which implies (\ref{P-cota}).
\end{proof}
\begin{lemma}\label{exp-bounds1}
There exist positive constants $\lambda$ and $D$ such that 
\begin{equation}
    \|T(t,s)P(s)\|\leq D e^{-\lambda(t-s)} \mbox{ for } t\geq s.\nonumber
\end{equation}
\end{lemma}
\begin{proof}
Let $\phi\in E(s)$ and $u_t=T(t,s)\phi$ for $t\geq s$.  We define the function $x:\mathbb{R}\to\mathbb{R}^n$ as
\begin{equation}
    x(t)=\left\{\begin{array}{lr}
        u(t)-\psi(t) &  \mbox{if } t\geq s-r, \\
        0 &  \mbox{if } t<s-r
    \end{array}\right..\nonumber
\end{equation}
Clearly $x$ is bounded and continuous on $\R$ , hence $x\in L^p_\mu$.  
Moreover:
\begin{enumerate}
\item  For $t\geq s$, $x$ is differentiable and  
\begin{eqnarray}
    x'(t) &=& u'(t) - \psi'(t)\nonumber\\
        &=&  L(t)u_{t}\nonumber\\
        &=&L(t)(x_t+\psi_t)\nonumber\\
        &=&L(t)x_t + g(t),\nonumber
\end{eqnarray}
\item For $t\leq s$, $x(t)=0$, and as $g(t)=0$ in the same interval, we have also that $x'(t)=L(t)x_t + g(t)$.
\end{enumerate}
In conclusion $\Lambda_L x=g$ and then $x\in W^{1,p}_{\mu}$. This implies that, in particular for $t\geq s$,  
\begin{eqnarray}
   \|u_t\| &=&\|x_t+\psi_t\| \nonumber\\
    & \leq & \sup\limits_{t\in\mathbb{R}}|x(t)| + \|\phi\|\nonumber\\
    &\leq & \|\Lambda_L^{-1}\|\beta\|\phi\| + \|\phi\|\nonumber\\
    &=&\gamma_0 \|\phi\|.\label{u-cota1}
\end{eqnarray}
Note also that if $\tau\in [s,t]$, repeating the same argument we can write 
\begin{equation}\label{u-cota1.1}
\|u_t\|\leq \gamma_0\|u_\tau\|.
\end{equation}
Now we shall prove that there exists $l\in\mathbb{N}$ (independent of $s$ and $\phi$) such that
\begin{equation}\label{u-cota2}
    \|u_t\| \leq \frac{1}{2}\|\phi\| \qquad \mbox{ for } t-s\geq l.
\end{equation}
If there is some  $t_*>s$ with satisfying 
\begin{equation}\label{medio}
\|u_{t_*}\|>\frac{1}{2}\|\phi\|,
\end{equation}
then we can find  $t_0\in[t_*-r,t_*]$ satisfying

\begin{equation}\label{u-cota0}|u(t_0)|=\|u_{t_*}\|>\di\frac{\|\phi\|}{2}
 \end{equation} 
and  in consequence \begin{equation}\label{u-cota0.1}\|u_{t_0}\|\geq|u(t_0)|>\di\frac{\|\phi\|}{2}
 \end{equation}

By other hand, from (\ref{u-cota1}), (\ref{u-cota1.1}) and (\ref{u-cota0.1}), for $s\leq \tau\leq t_0$, we obtain
\begin{equation}\label{u-cota3}
    \frac{1}{2\gamma_0}\|\phi\| < \|u_{\tau}\| \leq \gamma_0\|\phi\|.
\end{equation}

 Consider
\begin{eqnarray}
   I(t)&:=&\int_{-\infty}^{t}\chi_{[s,t_0]}(\tau)\|u_{\tau}\|^{-1}d\tau,\nonumber\\
    w(t)&:=&  u(t)I(t),\nonumber
\end{eqnarray}

Hence $w$ is continuous and
\begin{eqnarray}
    |w(t)|&\leq &|u(t)|\int_{s}^{t_0}\|u_{\tau}\|^{-1}d\tau\nonumber\\
          &\leq &\gamma_0\|\phi\|\int_s^{t_0} \|u_{\tau}\|^{-1}d\tau\nonumber\\
           &\leq& 2\gamma_0^2 (t_0-s),\nonumber
\end{eqnarray}
 this implies $w\in L^p_\mu$. Also
\begin{equation}
w'(t)=L(t)w_t+h(t),
\end{equation}
where 
\begin{equation}
    h(t) = I(t)L(t)u_t-L(t)(u_tI_t)+u(t)\chi_{[s,t_0]}(t)\|u_{t}\|^{-1},\nonumber
\end{equation}
for $t\in\mathbb{R}$. If $t\leq s$, $h(t)=0$ and if $t\geq s $ is easy to see that $h$ is bounded and piecewise continuous, in fact using  (\ref{u-cota3}) we have 
\begin{eqnarray}\label{h-cota}
    |h(t)|&=&|L(t)[I(t)u_t-u_tI_t]+u(t)\chi_{[s,t_0]}(t)\|u_{t}\|^{-1}| \nonumber\\
       &\leq & \beta \sup\limits_{-r\leq \theta\leq 0} |I(t)u(t+\theta)-u(t+\theta)I(t+\theta)|+1\nonumber\\
       &\leq & \beta \sup\limits_{-r\leq \theta\leq 0} \left\{ |u(t+\theta)|\int_{t+\theta}^t \|u_\tau\|^{-1} d\tau\right\}+1\nonumber\\
              &\leq & \beta \sup\limits_{-r\leq \theta\leq 0} \left\{\int_{t+\theta}^t \frac{u(t+\theta)}{\sup\limits_{-r\leq z\leq 0}u(\tau+z)} d\tau\right\}+1\nonumber\\
               &\leq & \beta \sup\limits_{-r\leq \theta\leq 0} \left\{\int_{t+\theta}^t d\tau\right\}+1\nonumber\\
               &=&\beta r+1.\nonumber
\end{eqnarray}
Hence $h\in L^p_\mu$ and $w\in W^{1,p}_\mu$. It allows to find the following bound independent of $t_0$ and $s$
\begin{equation}
\sup\limits_{t\in\R}|w(t)|\leq \|\Lambda^{-1}_L\|(\beta r+1)
\end{equation}

From (\ref{u-cota0}) and  (\ref{u-cota3}) we obtain
\begin{eqnarray}
    \|\Lambda_L^{-1}\|(\beta r+1) &\geq & |w(t_0)|\nonumber\\
    &\geq & |u(t_0)|\int_s^{t_0}\|u_{\tau}\|^{-1}d\tau\nonumber\\
    &\geq & \frac{\|\phi\|}{2}\frac{(t_0-s)}{\gamma_0\|\phi\|}\nonumber\\
    &\geq & \frac{(t_0-s)}{2\gamma_0}.\nonumber
\end{eqnarray}
Hence 
\begin{eqnarray}
t_*-s&=&t_*-t_0+t_0-s\\
    &\leq &2\gamma_0\|\Lambda_L^{-1}\|(\beta r+1)+r.\nonumber
\end{eqnarray}
Taking
\begin{equation}\label{l-cota}
    l>2\gamma_0\|\Lambda_L^{-1}\|(\beta r+1)+r
\end{equation}
the property (\ref{u-cota2}) holds.

Now, for any  $\phi\in C$  and    $t\geq s$, we write $t-s=kl+\tau$, with $k\in\mathbb{N}\cup \{0\}$ and $0\leq\tau<l$. Appliying (\ref{u-cota1}) and (\ref{u-cota2}) we have

\begin{eqnarray}
    \|T(t,s)P(s)\phi\| &=& \|T(s+kl+\tau,s)P(s)\phi\|\nonumber\\
    &\leq &\frac{1}{2^k}\|T(s+\tau,s)P(s)\phi\|\nonumber\\
    &\leq &\frac{\gamma_0}{2^k}\|P(s)\phi\|\nonumber\\
    &\leq &\gamma_0^2 2^{-k}\|\phi\|\nonumber\\
    &\leq &2\gamma_0^2 2^{-(t-s)/l}\|\phi\|\nonumber\\
   &=& 2 \gamma_0^2 e^{-\lambda(t-s)}\|\phi\|,\nonumber
\end{eqnarray}

with $\lambda=\di\frac{\ln(2)}{l}$, and  the proof is completed.
\end{proof}

\begin{lemma}\label{exp-bounds2}
There exist positive constants $\lambda$ and $D$ such that 
\begin{equation}
    \|\overline{T}(t,s)Q(s)\| \leq D e^{-\lambda(s-t)} \mbox{ for } t\leq s.\nonumber
\end{equation}
\end{lemma}
\begin{proof}
The proof is analogous to the previous one, defining
for $s\in\mathbb{R}$ and $\phi\in F(s)$ the functions
\begin{equation}\label{psi-function}
\psi(t)=\left\{\begin{array}{ll}\phi(t-s),&\mbox{ if } t\in [s-r,s],\\
                                                                   \phi(-r),&\mbox{ if } t\in (-\infty,s-r)
                                                  \end{array}\right.
                                                  \end{equation}
                  and $g:\mathbb{R}\to\mathbb{R}^n$ by 
\begin{equation}
  g(t)=\left\{
    \begin{array}{lr}
        0 & \mbox{ if } t> s \\
        L(t)\psi_t & \mbox{ if } t\leq s 
    \end{array}\right..\nonumber  
\end{equation}
\end{proof}
Combining the previous lemmas, we conclude that the equation (\ref{eq0}) has an exponential dichotomy on $\mathbb{R}$.

\end{document}